\documentclass[12 pt]{article}        	
\usepackage{amsfonts, amssymb}
\usepackage{amsmath}
\usepackage{graphicx}
\usepackage{float}
\usepackage{esint}
\usepackage[left=2cm,right=2cm,top=2cm,bottom=2cm]{geometry}
\newcommand{\real}{\mathbb R}

\newcommand{\n}{\mathbb N}

\usepackage{units}

\title{Modeling the Impact of NATO Policy Actions on Taliban Disinformation Campaigns with Lotka-Volterra Models}

\date{March 2025}

\begin{document}

\maketitle
\begin{center}
    \textbf{Author}

    Timothy Tarter, James Madison University

    Department of Mathematics

    www.tartermathematics.com

    $\newline$

    Bella Santos, James Madison University

    Department of Public Policy    
\end{center}

\begin{abstract}
    In this paper, we use an adaptation of the Lotka-Volterra Predator-Prey framework (via Affili Et. Al.) to model the relationship between Taliban disinformation campaigns and NATO's attempts to counter them. We find that aggression is a significantly preferred strategy over stagnation; this result should prompt decisionmakers to take strong stances against disinformation campaigns while maintaining a healthy defense against retribution. This leads us to propose the NATO Afghan Youth for Change program, which provides scholarships to individuals throughout the Afghan diaspora to become NATO influencers - a role which will both advance their personal goals and mitigate the spread of disinformation in what may otherwise be communities targeted by the Taliban's disinformation campaigns.
\end{abstract}

\section{Overview}
This research came about through a NATO-sponsored hackathon at William and Mary centered around countering disinformation threats. Our stream focused on the Taliban's attempts to percolate radical ideology online and expand terrorist cells in NATO countries through disinformation. Our initial work focused on perturbation theory of a simple dynamical system, with the understanding that the true system was significantly more complex than our preliminary model. This subsequent work builds on our previous framework by delving into the intricacies which we were previously unable to address due to the time constraint of the hackathon.

\subsection{Problem Statement}
In the modern era of warfare, we face a set of multifaceted threats that have no historical precedent for mitigation. From payload drones to cyber threats to radical disinformation campaigns, creating an effective mitigation strategy for a threat is more than just problem-solving, it is situational awareness, intelligence analysis, and cognizance of implications. Narrowing in on our threat space, online disinformation campaigns from radical groups such as the Taliban are becoming more common for three reasons. First, disinformation campaigns are cost-efficient because they only require an internet connection and a VPN. Next, disinformation campaigns are unlikely to cause unexpected blow-back due to the isolated yet expansive nature of the Taliban's operations. Finally, disinformation campaigns don't require radical recipients to propagate the threat - lack of education, gullible readers, and ``shock factor" pass disinformation through the western populous like a virus through a crowd. Because of this tendency for threat replication, we want to provide a quantitative analysis to support certain policy strategies which will (provably) mitigate this threat in the three-to-five years.

\subsubsection{Defining Terminology for the Threat Landscape}
To accurately address the threat landscape, we must first clarify some common definitions, which are taken \textit{verbatim} from NATO's most recent policy report on disinformation.$^{[1]}$
\begin{itemize}
    \item Misinformation is therefore excluded from the definition of information threats. This phenomenon describes false or inaccurate information spread without malicious intent. However, the effects of misinformation can still be harmful, making proactive communications and easy access to facts even more important.
    \item Disinformation – false or inaccurate information spread deliberately to manipulate the opinions and actions of others.
    \item Propaganda – information designed to manipulate a specific target audience toward a particular behaviour or belief, often as part of a prolonged campaign by a state actor with a political agenda.
    \item Information Manipulation and Interference by Foreign Actors (IMIF) – a pattern of behaviour that threatens or has the potential to negatively impact values, procedures and political processes in a target country. Such activity is mostly non-illegal, but is manipulative in character, conducted in an intentional and coordinated manner by state or non-state actors (including their proxies inside and outside of their own territory).
    \item Hybrid warfare – use of military and non-military as well as covert and overt means (including disinformation, cyber attacks, economic pressure, deployment of irregular armed groups, and use of regular forces) to blur the lines between war and peace, sow doubt in the minds of target populations, and destabilise and undermine societies.
\end{itemize}

\section{Strategy Selection}
Our quantitative modeling will revolve around strategy selection, specifically, examining how different policy strategies affect: 
\begin{itemize}
    \item the degree of `eradication' of Taliban disinformation campaigns by a specific strategy
    \item the rate at which that `eradication' occurs given that strategy
    \item how combining strategies (on both sides) affects the timeline for threat development. 
\end{itemize}
We will characterize these strategies qualitatively here in order to both justify and expedite our model selection process for generating simultaneous equations in section [4].

\subsection{Current NATO Strategies}
According to a policy report$^{[1]}$ from December 2024, NATO currently employs short, medium, and long-term strategies based on five key functions: 
\begin{enumerate}
    \item influencer campaigns 
    \item understanding the information environment
    \item preventing the effectiveness of information threats
    \item containing and mitigating specific information incidents
    \item recovering stronger by learning lessons from information threats
\end{enumerate}

$\newline$
However, while these tenants comprise NATO’s overall approach to disinformation, it is important to note that these are generalized information strategies. NATO’s strategy pertaining to counter-terrorism is largely centered around physical military capabilities and hybrid warfare, not necessarily the propagation of false information that incentivizes recruitment or instills positive perceptions of groups like the Taliban. Nonetheless, it is important to identify current NATO strategies to combat disinformation to see how anti-terrorism efforts can fit into existing frameworks. 

\subsubsection{Influencer Campaigns}
In 2022, NATO founded the `Protect the Future' campaign,$^{[2]}$ an initiative which allowed NATO to inform younger audiences about events in international relations and public policy through popular social media ``influencers." Since the program's inauguration, influencers have attended NATO Summits, visited aircraft carriers, hosted gaming tournaments, and even launched a graphic novel - all for the purpose of greater engagement with groups who may otherwise fall prey to disinformation campaigns. In the program's purpose statement, NATO details that this is accomplished by, ``giving young people a voice in telling NATO's story in their own way."$^{[2]}$ 

$\newline$On a less idealistic note, the program is a tractable method for approaching audiences which are easily influenced. According to the Pew Research Center$^{[3]}$ in 2024, roughly a quarter of Americans said they often get their news through social media. 
\begin{figure}[H]
    \centering
    \includegraphics[width=0.4\linewidth]{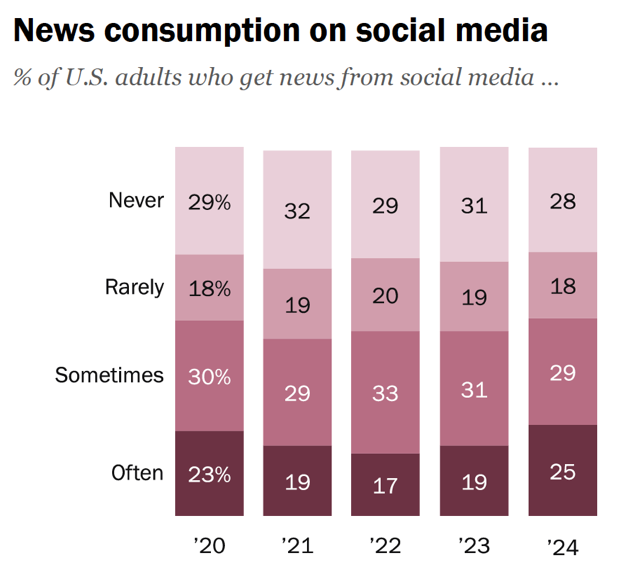}
    \caption{News \& Media Consumption (Pew, 2024)$^{[3]}$}
    \label{fig:enter-label}
\end{figure}

$\newline$
Moreover, according to their study, for users of Facebook, X (formerly Twitter), and TikTok, between 52\% and 59\% of users regularly got their news from the respective social media websites in 2024. 

\begin{figure}[H]
    \centering
    \includegraphics[width=0.75\linewidth]{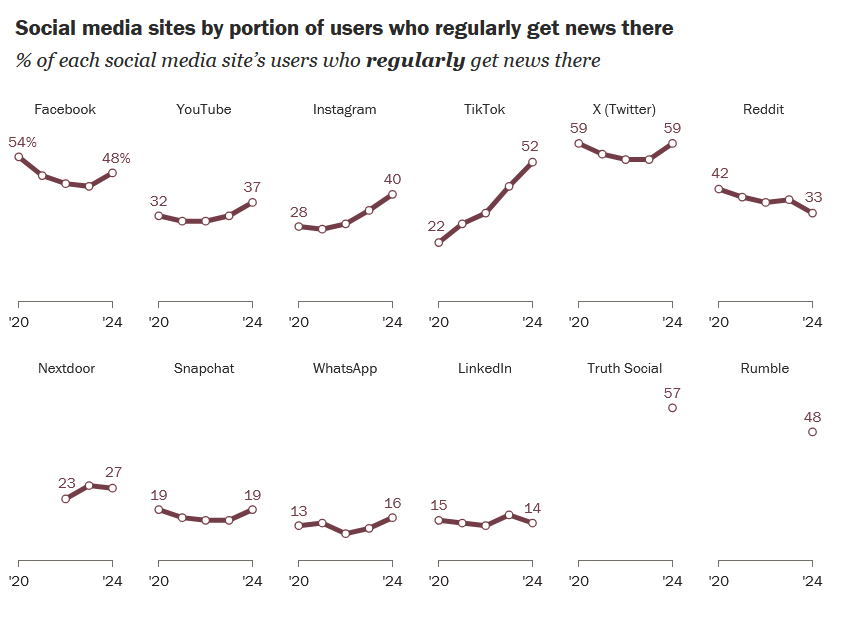}
    \caption{Regular `News' Consumption (Pew, 2024)$^{[3]}$}
    \label{fig:enter-label}
\end{figure}

$\newline$
In light of these statistics, it is clear that NATO funded influencer campaigns are a robust way to counter disinformation by directly spreading correct information to the populous.   

\subsubsection{Understanding the Information Environment}
In addition to any other strategy, for NATO to effectively combat the spread of disinformation, the Alliance must accurately understand the information environment and active narratives that threaten its security. NATO internally, as well as externally and in collaboration with members’ own intelligence services, tracks and analyzes relevant information. Central to this is NATO’s Information Environment Assessment (IEA) capability,$^{[4]}$ which aims to track both friendly and adversarial information. In 2018, NATO’s Operational Experimentation (OpEx) conducted a series of experiments to develop NATO’s IEA capabilities. Their goal was to reach initial operating capability in June 2019, and final operating capability in 2022. While they did reach their initial operating capability goal in 2019, IEA has not yet reached final operating capability, and though it is in use it continues to undergo development and refinement.

$\newline$
The IEA strategy identifies that it is not only disinformation that is critical to the “understand” function; rather, they aim to critically examine “own, earned, and hostile” information. The Alliance seeks to analyze how the information they put out, news coverage they receive, and information by malign actors are respectively received and proliferated. As of 2023, NATO’s Strategic Warfare Development Command highlighted two aspects of IEA capabilities: technology and personnel.$^{[5]}$ There is a significant technological component: utilizing AI and machine learning techniques allows NATO to continuously analyze the information environment and identify groups or behaviors that pose a threat. Using this technology now will support automated predictive analytics in the future, allowing the Alliance to optimize decision-making at even faster rates. However, there is also a need for qualified personnel with backgrounds in data analysis, especially for those with military experience. NATO published that, “To achieve Full Operational Capability, a total of 48 new civilian positions and 18 reclassified military personnel positions must be allocated to Information Environment Assessment throughout the NATO Command Structure and NATO Headquarters Public Diplomacy Division to build the necessary teams to achieve Strategic Communications objectives.”$^{[5]}$ This addition would be further complemented by continued development of IEA doctrine and policy.  

\subsubsection{Preventing Effectiveness}
The bulk of NATO’s campaign to counter disinformation falls under the category of prevention, specifically, inoculation against false information and spreading positive perceptions of the Alliance. 

$\newline$
The primary aim of this strategy is to introduce credible information preemptively to prime the public to critically intake disinformation published by malign actors. This occurs through proactively sharing information with the public through coordination between allies. NATO’s Rapid Response Group (NRRG) provides early warning of informational threats and coordinates information sharing to respond. The NRRG complements other intelligence sharing forums that NATO is a party to. For instance, NATO is an observer to the G7 Rapid Response Mechanism. The G7 RRM is a collective response mechanism that allows G7 members and its observers to develop shared databases to track foreign threats, identify collective responses to disinformation and security threats, and communicate its findings to the public. In the G7 RRM’s 2022 Annual Report, they acknowledge NATO’s preventative work to counter disinformation. For example in 2021, the NATO Alliance declassified intelligence regarding Russia’s military buildup and intent to stage a full-fledged invasion of Ukraine, even warning against potential false flag operations.$^{[6]}$ This intelligence was shared with the public through coordination with allies to increase resilience in public audiences against Russian disinformation. Public statements coming directly from figures like NATO’s Secretary General and High Ranking Officials within alliance states gave this “pre-bunking” an added layer of legitimacy and so the campaign was largely effective.

\subsubsection{Containing and Mitigating Incidents}
Some intelligence professionals argue that engaging with false facts to disavow them only contributes to their proliferation by giving them a platform and an air of legitimacy. Thus, NATO primarily sticks to rules of non-engagement with such narratives, however, the Alliance does attempt to “debunk” certain false information once it becomes widespread and a clear harm to the public. Methods regarding containing and mitigating the proliferation of disinformation commonly include public communication via official statements to the media and posts on the NATO website and social media accounts.$^{[1]}$

$\newline$
NATO’s debunking capabilities are another aspect of their approach that was especially highlighted following the Russian invasion of Ukraine. In 2014, NATO established a page on their website titled, “Setting the Record Straight.” Information contained there was primarily in response to Russian disinformation regarding the illegal annexation of Crimea. $^{[8]}$ However, following the invasion of Ukraine, the website was updated to respond to prevailing Russian claims justifying their military aggression. Each of these claims and responses that NATO has made follow a similar pattern: calling out inaccuracies in Russian claims and emphasizing their hypocrisy. 

\subsubsection{Recover Stronger}
The final pillar of NATO’s agenda to counter disinformation is the shortest and least concrete. It does not necessarily stand on its own, but rather exhibits an accumulation of learned experiences from the strategies outlined above. In many ways, there is substantial overlap between the first approach, “understanding the information environment,” and the ability to recover stronger. As mentioned, one of NATO’s goals in establishing their IEA capability is to develop future automated predictive analytics that allow the Alliance to more quickly understand evolving narratives that threaten NATO’s security. Remaining cognizant of the manipulative tactics, techniques, and procedures (TTPs) through which malignant actors seek to sow distrust or discontent means that NATO is able to more effectively respond to future attacks of all types. 

$\newline$
NATO’s collaboration with allies and other international organizations also contributes to their ability to recover stronger by creating a supportive network to collectively respond to amorphous threats. Through partnerships with entities like the G7 RRM, NATO participates in forums that actively work to learn from previous intelligence events while coordinating group action rather than leaving individual member states to go it alone.

\subsection{Current Taliban Strategies}
After they began to lose power in the 1990s, the Taliban turned to social media and the internet as a cost-effective method of spreading desired narratives. Social media use was particularly important to the resurgence of the Taliban following the U.S. withdrawal from the country in 2021. During their fight to take back Kabul, they regularly posted to social media highlighting their territorial gains.$^{[9]}$ This influenced the Afghan public and furthered the narrative that the Taliban was sustaining control of large portions of the country. In many ways, the use of social media allowed the group to overstate their successes and portray themselves as a swift and decisive military force. Additionally, they posted videos of Afghan government soldiers defecting to join them, or surrendering to their troops. This decreased public confidence in the Afghan army and reduced morale amongst the remaining soldiers. Ultimately, these disinformation campaigns heightened public fear regarding another Taliban regime and contributed to the rapid deterioration of the American-backed Afghan government. 

$\newline$
Social media provides a low-cost method of spreading disinformation within Afghanistan, as well as to the broader world. The Taliban frequently employs narratives referring to themselves as the rightful rulers of Afghanistan. During the years of U.S. presence in the country, they often framed an existential battle between the Afghan government (a puppet of the West as they claimed), aided by the Americans, versus their genuine understanding of Islam and Afghan culture. For Afghan youth who have lived through years of conflict and domestic turmoil, these narratives can be particularly appealing. Young people between the ages of 15-24 who spend large amounts of time online are especially susceptible to radicalization, and terrorist groups frequently employ social media as a recruitment tool for Muslim youth.$^{[9]}$

$\newline$
Social media is also advantageous because information provides an accessible way to communicate with audiences globally. The Taliban has used highly personalized social media messages in order to advance their cause. For instance, before the current iteration of the Taliban rose to power, former Taliban spokesperson Zabiullah Mujahid directly responded to various comments left on their social media accounts, clarifying the Taliban’s perspectives on given topics.$^{[9]}$ Social media audiences come from all corners of the world, and the Taliban has sought to exploit this opportunity. Taliban social media accounts have published content in multiple languages, including English, Urdu, Turkish, Pashto, and Persian - a clear acknowledgment of their desire for extensive and expansive reach.$^{[9]}$

\subsection{Proposed NATO Strategies}
In this paper, we propose the NATO Afghan Youth for Change (AYFC) program, which would grant individuals throughout the Afghan diaspora scholarships in exchange for being NATO influencers to their communities. Beyond the larger scope of social media, the AYFC hopes to provide validated information to groups which the Taliban may otherwise seek to convert to its radical agenda. While this may seem like a very ``aggressive" strategy on the part of NATO, due to the rise in extremism throughout Europe and other parts of the Afghan diaspora, we hold that an aggressive strategy is necessary to limit the spread of terrorist ideology. Beyond just this qualitative theory, however, we further substantiate the merit of such a program in chapters five and six with simulation models.

\section{Literature Review}
\subsection{A Brief Overview of Differential Equations}
As this is a paper in an international relations publication (not a math journal), we feel that we should include a short primer for policy analysts on what differential equations are so that they can understand the claims we make in future sections. While we are going to go into full rigorous detail of the models employed in these future sections, we invite those without a math background to still read this section of our publication. Differential equations are a surprisingly intuitive corner of mathematics - specifically, they are equations relating rates of change. We all encounter them throughout our lives, from observing the physical world around us, to managing money, to even controlling our stress levels. For example, consider the rate of change of one's stress. It may be affected by any number of variables, but let's say for our example that it is only dependent on the current amount of stress we have, the amount of money in our bank account, and the amount satisfaction we have from our homes, jobs, and hobbies. This is an example of a differential equation. We haven't written it nearly as formally as we will for our other ODEs (Ordinary Differential Equations) in this paper, but it doesn't change the fact that it is one. Often times though, we don't just want to know what the rate of change of something at a certain point in time is, we want to know what the actual quantity of the thing is at that point in time. This is where we get into \textit{solving} ODEs - a significantly more complicated task if we only accept symbolic solutions. 

$\newline$We are used to the idea from highschool and college that equations should be solvable for exact numbers, but this often isn't the case - especially for complicated ODEs and PDEs (partial differential equations). When this happens, we typically say that a system has ``no closed-form (analytic) solution" and employ both qualitative methods and quantitative methods of studying the system. One such qualitative method is examining the ``phase-plane" (which is just examining a plot of the state vector components against one another across time); this is very useful for our specific application because we can study cyclic phenomena and see clearly when a system has been perturbed from a cycle. Many quantitative methods are centered around approximation (such as the Parker-Sochacki method), as well as mathematically proving what can and can't happen. We'll see these show up when we discuss the dynamics of the system in section (5). 

$\newline$Interpreting differential equations is often the easy part - after all, we have an intuition of what generally \textit{should} happen from real world experience. Moreover, since our model has been built on that intuition, we can usually pinpoint and prove intricacies about our system that we would otherwise miss if we just focused on things like mean response (not to say that econometric modeling isn't helpful, just a different tool). 

\subsection{Classic Lotka-Volterra Model}
The Lotka-Volterra model is a system of nonlinear differential equations that was developed to explicitly model predator-prey population dynamics for mathematical biology. While it is a very simple model, it has shown the test of time as a robust framework for modeling ecosystems. 

$\newline$The assumptions that the model make are as follows (GeeksForGeeks, 2024):$^{[11]}$
\begin{itemize}
    \item The population of prey increases either linearly or exponentially in the absence of predators.
    \item The availability of prey is the only factor limiting the predator's population expansion.
    \item There is no age structure, constant reproduction, or atomic homogeneity between predator and prey.
    \item The rate of encounters between predators and prey is directly correlated with the rate of predation.
    \item The predator has a density-independent, constant mortality rate.
\end{itemize}

\subsubsection{Parameters}
The model takes into account the following parameters:$^{[11]}$
\begin{itemize}
    \item $P$ = prey population
    \item $Q$ = predator population
    \item $a$ = birth rate of prey
    \item $b$ = predation rate (the rate at which the predators kill the prey)
    \item $\delta$ = efficiency rate (the rate at which eating prey creates new predators)
    \item $n$ = natural death rate of predators (the rate at which predators will die if there are no prey to eat)
\end{itemize}

\subsubsection{Model}
The formalization of the system of equations follows the rates of change of both the populations. For the prey model, we get:
\begin{equation}
    \frac{dP}{dt} = P(a-bQ)
\end{equation}

$\newline$
Likewise, for the predator model, we get:
\begin{equation}
    \frac{dQ}{dt} = Q(-n + \delta P)
\end{equation}

\subsubsection{Model Discussion}
When we initially began this research in the WMGIC x NATO 2024 Hackathon, this was the model that we employed to study the online threat landscape. While it served well enough for the short time period we were given to come up with a project, there were some qualities it lacked. For one, it assumes that the predators have no other food - so in a way, we expect the dynamic to either be a coexistence cycle or end in an extinction event. We were still able to study perturbations that led to extinction events, but we want a model that lets one population live after the other dies. Otherwise, the model will have a bias of omission. Moreover, the classical Lotka-Volterra model assumes homogeneity of populations. While the Taliban is largely composed of a homogeneous population of Afghani citizens, we may want a model that can be expanded to model online threats with a less homogeneous population. Finally, the largest assumption is homogeneity of strategy - which is something our paper deals quite heavily with. Thus, this model isn't quite suitable for our purposes. That said, these needs motivate our next model.

\subsection{Affili Et. Al.'s Civil War Model}
In 2024, Affili Et. Al. published a monograph titled, ``A New Lotka-Volterra Model of Competition with Strategic Aggression: Civil Wars When Strategy Comes Into Play"$^{[16]}$ under Springer, a book which made significant developments on systems with shared resources - such as publicity and public belief (for our use case). The model doesn't deviate too heavily from the classic Lotka-Volterra model, but does introduce the novel idea of a variational ``strategy" parameter. They do an extremely thorough job of proving the extensive dynamics of their model, and we will reuse many of their results, theorems, and lemmas. We ultimately settle on this model because of strong similarities between the assumptions of Affili Et. Al., and our own. 

\subsubsection{Parameters}
The factors that the `Civil War' model takes into account are:
\begin{itemize}
    \item a = the portion of the population fighting = level of aggression = strategy
    \item p = the second group's fitness to resources 
    \item $c_i$ = $\frac{endured}{inflicted}$ kills for population `i' (note: $c_1 = \frac{1}{c_2}$)
    \item u = the population density of the first population = $\frac{pop_1}{pop_1 + pop_2}$
    \item v = the population density of the second population = $\frac{pop_2}{pop_1+pop_2}$
\end{itemize}

\subsubsection{Model}
The system of differential equations is as follows:
\begin{equation}
    \frac{du}{dt} = u(1-u-v) - acu
\end{equation}
\begin{equation}
    \frac{dv}{dt} = pv(1-u-v) - au
\end{equation}

\subsubsection{Model Discussion \& Important Results}
A big element worth noting about this model is that everything is measured in density - so the phase plane is constrained to $[0,1]\times [0,1] \subseteq \real^2$. A great deal of the monograph is devoted to partitioning this solution space into finite and non-finite time convergence solutions - namely, for any smooth curve $\mu$ that enters and exits $G = [0,1]\times [0,1] \subseteq \real^2$ exactly once, we can partition G into two sets, $\beta$ (the set of infinite time solutions) and $\varepsilon$ (the set of finite time solutions). Additionally, we note that $\mu$ is a stable solution.

\paragraph{G is Partitioned into Finite and Non-Finite Time Solutions by a Smooth Curve, $\mu$:}

$\newline$
It is proved in the monograph that certain trajectories cannot enter or exit through $\mu$ based on their starting position. Formally, let S = $\{(u,v) \in [0,1]\times [0,1] $ $|$ $v < g(u)\}$, where $g$:[0,1] $\to$ [0, $+\infty)$ is a monotone increasing, continuous function such that there exists a finite (possibly empty) set Z $\subset$ [0,1] for which g $\in$ $C^1$([0,1]$\setminus$Z), and in addition, the limits
\begin{equation*}
    \lim_{u \to z^\pm}g'(u)
\end{equation*}
exist for all $z \in Z$ (possibly distinct and possibly equal to $+\infty$). Let N($\hat{u}, g(\hat{u})$) be the set of outward unit normal vectors to $\partial$S at ($\hat{u}, g(\hat{u})$), defined by
\begin{equation}
    v = (-\frac{g'(\hat{u})}{\sqrt{1+(g'(\hat{u}))^2}}, \frac{1}{\sqrt{1+(g'(\hat{u}))^2}}).
\end{equation}
Let $\{I_k\}_{k\in A}$, where A is the set of admissible strategies, and $\{k_j\}_{j \in C}$, where C is a finite subset of $\n$, be two finite collections of disjoint intervals in $[0, 1]$, the $I_k$ being open and the $K_j$ being closed (possibly coinciding with singletons) in the topology induced by $[0,1]$, such that
\begin{equation}
    \bigcup_{k\in A}I_k \cup \bigcup_{j \in C}K_j = [0,1].
\end{equation}

$\newline$
Suppose that for all $\hat{u} \in \bigcup_{k \in A}I_k$ it holds that
\begin{equation}
    \min_{v\in N(\hat{u}, g(\hat{u})))} (F(\hat{u}, g(\hat{u}), G(\hat{u}, g(\hat{u})))\cdot v > 0
\end{equation}

$\newline$
and that for all $\hat{u} \in \bigcup_{j \in C}k_j$ it holds that 
\begin{equation}
    \min_{v \in N(\hat{u}, g(\hat{u}))} (F(\hat{u}, g(\hat{u})), G(\hat{u}, g(\hat{u})))\cdot v = 0.
\end{equation}

$\newline$
Then no trajectory enters S through a point of the set
\begin{equation}
    \partial S \cap \{(u,v) : u \in [0,1], v = g(u)\}.
\end{equation}

$\newline$
The proof of this is in \textbf{no way trivial}, and can be found on pages 32 - 36 of the monograph. That said, this proof is extremely important - this implies that non-finite time solutions never become finite time solutions. 

$\newline$
Similarly, it is shown (in a much shorter proof no less) that no finite-time solutions ever become non-finite time solutions, and thus, we have a complete partition of G = $[0,1] \times [0,1]$.

\paragraph{G is Partitioned into Four Basins of Attraction, and No Trajectory Exits $A_3$ if ac $\not=$ 1:}

$\newline$
It is also proved that there exist four unique basins of attraction, $A_1, A_2, A_3,$ and $A_4$, which are divided by the curves 
\begin{equation}
    \frac{dv}{dt} = 0
\end{equation}
which is given by
\begin{equation}
    \sigma(v) = 1 - \frac{pv^2+a}{pv+a}
\end{equation}
and
\begin{equation}
    \frac{du}{dt} = 0
\end{equation}

$\newline$
The formalization of these partitions are as follows: 
\begin{itemize}
    \item $A_1 = \{(u,v) \in [0, 1] \times [0,1] $ $|$ $\frac{du}{dt}\leq 0, \frac{dv}{dt} \geq 0\}$
    \item $A_2 = \{(u,v) \in [0, 1] \times [0,1] $ $|$ $\frac{du}{dt}\leq 0, \frac{dv}{dt} \leq 0\}$
    \item $A_3 = \{(u,v) \in [0, 1] \times [0,1] $ $|$ $\frac{du}{dt}\geq 0, \frac{dv}{dt} \leq 0\}$
    \item $A_4 = \{(u,v) \in [0, 1] \times [0,1] $ $|$ $\frac{du}{dt}\geq 0, \frac{dv}{dt} \geq 0\}$
\end{itemize}

$\newline$
It is then proved that if ac $\not= 1$, that no trajectory exits $A_3$ - a very fruitful and interesting proof, especially for analyzing systems efficiently.

\paragraph{If ac = 1 and the Initial Condition Lies on $\mu$, its $\omega$-Limit is \{$(0,0)$\}:}
$\newline$
Affili Et. Al. defines the $\omega$-limit, given $(u_0, v_0) \in [0,1]\times[0,1]$, as:
\begin{center}
    $\omega(u_0, v_0) = \{(x,y) \in \real^2 : \phi_{(u_0,v_0)}(t) \in [0,1]\times[0,1],$
    
    $\newline$
    $\textrm{for all } t \geq 0,$

    $\newline$
    $\textrm{ and there exists }\{t_i\}_{i\in \n} \textrm{ with } t_i \to -\infty$,

    $\newline$
    $\textrm{and} \lim_{i \to +\infty}\phi_{(u_0, v_0)}(t_i) = (x,y)\}$
\end{center}

$\newline$In short, if ac = 1 and the initial condition lies on $\mu$, then the system will go to (0,0) in its limit. (Note: $\phi_{(u_0, v_0)}(t)$ is the trajectory in the phase plane.)

\paragraph{Equilibria Based on Values of ac:}
If ac $\in$ $(0,1)$:
\begin{itemize}
    \item (0, 0) is a source
    \item (0, 1) is a sink
    \item ($\frac{1-ac}{1+pc}pc$, $\frac{1-ac}{1+pc}$) is a saddle
\end{itemize}

$\newline$
If ac $> 1$:
\begin{itemize}
    \item (0, 0) is a saddle
    \item (0, 1) is a sink
\end{itemize}

$\newline$
If ac = 1:
\begin{itemize}
    \item (0, 0) has a positive eigenvalue and a null eigenvalue
    \item (0, 1) is a sink
\end{itemize}

\paragraph{How Fitness to Resources (p) Affects Strategy:}
$\newline$
If p = 1:
\begin{itemize}
    \item Strategy doesn't affect the outcome, just the time spent at war.
\end{itemize}

$\newline$
If p $<$ 1 (the first force is more efficient):
\begin{itemize}
    \item If u + v $\leq$ 1, then the resource market isn't over-saturated, and the first population should pick a small $a$. This is the same as not being very aggressive.
    \item If u + v $>$ 1, the first population should pick a very large $a$. This is the same as being very aggressive.
\end{itemize}

$\newline$
If p $>$ 1 (the second force is more efficient):
\begin{itemize}
    \item If u + v $\leq$ 1, pick a very large value of $a$, as a moderate or small $a$ value will cause a coexistence cycle.
    \item If u + v $>$ 1, pick $a$ = 0 until the resource market is no longer over-saturated, then pick a very large value of $a$.
\end{itemize}

\paragraph{Discussion:}
$\newline$
We want to focus on studying certain strategies already in use, some of which are constant in aggression, some of which are non-constant in aggression. Since each of these strategies are currently in use, we will have to beware of collinearity bias when examining the system dynamics from the model we build. That said, we can (to a certain degree) also model the interplay of these strategies, and attempt to build multi-faceted attacks to model. This will be accomplished in section (4).

\section{Modeling Strategies with Simultaneous Equations}
In this section, we want to model each individual strategy to find trial values of $a$ for our system of differential equations (given by the Civil War model in section (3)), from which we can draw policy conclusions. 

\subsection{Influencer Campaigns}
We want to characterize the aggression for these strategies over time using the following variables:
\begin{itemize}
    \item NATO Engagement $N_t$
    \item Identified Disinformation Campaign Engagement $D_t$
\end{itemize}
We derive the density of aggression for NATO engagement, $a_t$, as:
\begin{equation}
    a_t = \frac{N_t}{N_t + D_t}.
\end{equation}
$\newline$
Note: $a_t$ $\in$ [0,1] - this is to maintain continuity with the framework designed by Affili Et. Al.

$\newline$Engagement metrics for NATO can be found using OSINT analytics tools for different social media sites. According to Phalanx$^{[10]}$, NATO has a 1.13\% engagement rate with 52.4\% of their audience being real people, 2\% of their audience being influencers, and 44\% of their audience not engaging. Their expected authentic engagement per post is 8.3k users. 74\% of NATO's followers are male while only 26\% are female - a surprising discrepancy. Examining figure 3, we also see that their percentage engagement has been decreasing over the last four years in a roughly linear way. 

\begin{figure}[H]
    \centering
    \includegraphics[width=0.4\linewidth]{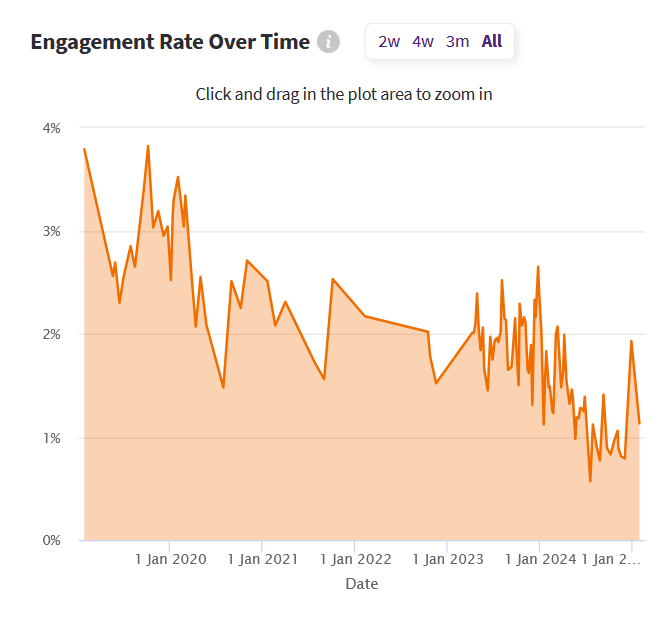}
    \caption{Engagement Time-Series Plot (Phalanx, 2025)}
    \label{fig:enter-label}
\end{figure}

$\newline$
In addition to the Phalanx report, our team created a database of NATO's Instagram likes over the last six months (July 2024 - January 2025), available online at www.tartermathematics.com. From this time series data, we find that there is a fairly drastic standard deviation of likes from post to post every month (see figure 4). Note: we chose months as our time-series lag unit as an arbitrary way to break up our sample. We also recognize the risk that omitted variable bias poses, especially given the width of our standard deviation.

\begin{figure}[H]
    \centering
    \includegraphics[width=.65\linewidth]{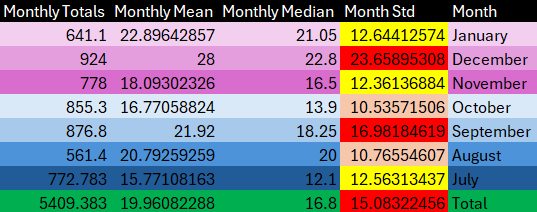}
    \caption{Four Number Summary by Month of NATO-Affiliated Instagram Likes (in Thousands)}
    \label{fig:enter-label}
\end{figure}

$\newline$Because of the magnitude of our standard deviation, we want to ensure that mean response is a viable measure of central tendency to parameterize $N_t$. If we fit a simple linear regression model (which is not our final method of parametrization but merely a diagnostic tool), we can show that our residuals are normally distributed and constant. Thus, $N_t$ = $\frac{\mu_i}{n}$ for $1 \leq i \leq n$ (over $n$ months of data) is an unbiased measure of central tendency. Let $\ell$ be the likes on NATO's Instagram, and $d$ be the date. Then our simple linear regression model will be of the form:

\begin{equation}
    \ell = \beta_0 + \beta_1d + \epsilon_i
\end{equation}

$\newline$Using the R output in figure 5 as justification, we use the parameters for our model:
\begin{equation}
    \ell = 151.425 + -0.75d
\end{equation}
\begin{figure}[H]
    \centering
    \includegraphics[width=0.75\linewidth]{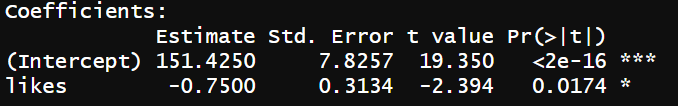}
    \caption{Parameter Significance}
    \label{fig:enter-label}
\end{figure}
$\newline$
To ensure that these significance tests are valid, we want to test for heteroskedasticity (non-constancy of residuals) using a Bresuch-Pagan test. We can quickly show that there is not heteroskedasticity in our model: 
\begin{figure}[H]
    \centering
    \includegraphics[width=0.55\linewidth]{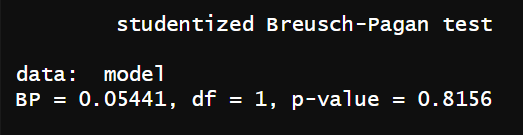}
    \caption{Breusch-Pagan Test for Heteroskdasticity}
    \label{fig:enter-label}
\end{figure}

$\newline$Since the null hypothesis is homoskedasticity, and since our p-value is \textit{well} above 0.05, we fail to reject the hypothesis of homoskedasticity and conclude that our residuals are constant. Thus, we know that our student t-tests for parameter significance are valid and our coefficients are efficient. Now we want to perform a Kolmogorov-Smirnoff goodness-of-fit test for normality of residuals as a means of deriving our empirical distribution that will yield the confidence interval for our mean response. Letting our null hypothesis be that our residuals are normally distributed and letting our alternative hypothesis be that they are not normally distributed, we can perform the K-S test in R:

\begin{figure}[H]
    \centering
    \includegraphics[width=0.6\linewidth]{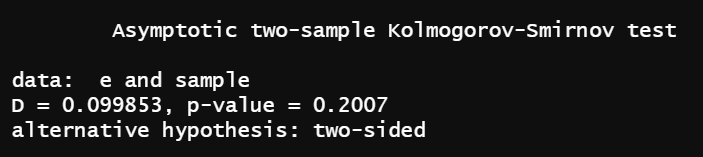}
    \caption{Kolmogorov-Smirnoff Goodness-of-Fit Test}
    \label{fig:enter-label}
\end{figure}

$\newline$Since p is greater than 0.05, we fail to reject our null hypothesis and conclude that our residuals are normally distributed. Visually, we can see this in figure 8 (empirical residuals are blue, simulated normal distribution is red).
\begin{figure}[H]
    \centering
    \includegraphics[width=0.65\linewidth]{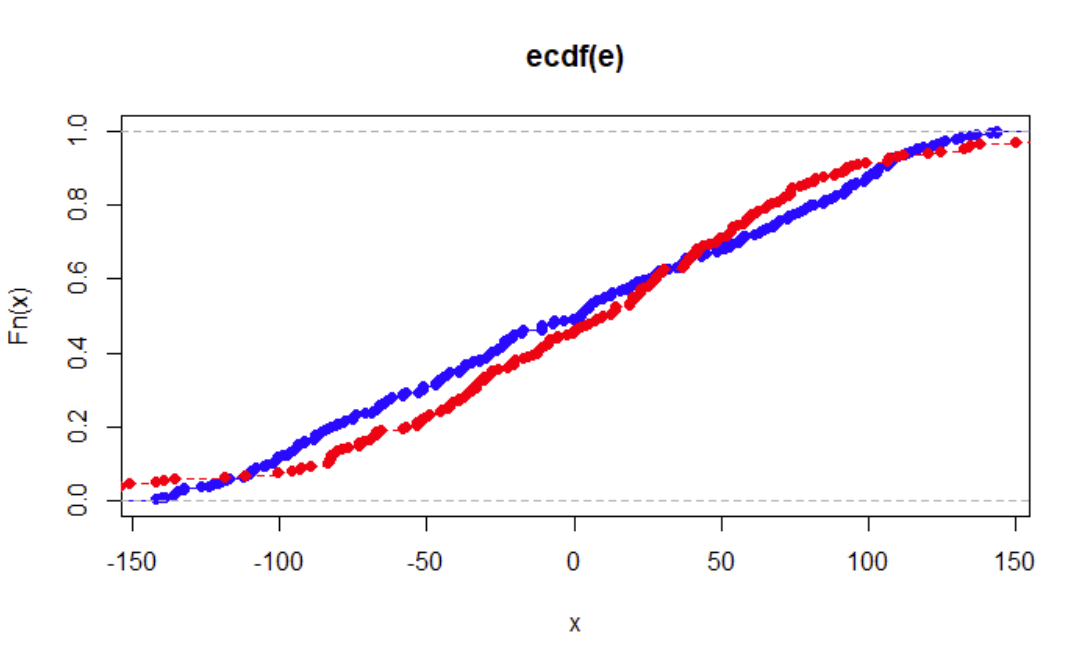}
    \caption{Goodness-of-Fit Plot}
    \label{fig:enter-label}
\end{figure}

$\newline$Since the residuals of our simple linear model are normally distributed and constant, and since our parameters are significant, we can reasonably use a normal distribution to model the change in engagement over time and induce a mean response function. 

$\newline$The big picture here is that even though our standard deviation is large in our time-series database, we can still model $N_t$ from \textbf{(13)} with the monthly mean from the data and the normal distribution.

$\newline$Next, we need to find a way to model the public's engagement with identified Taliban disinformation campaigns, $D_t$. According to researchers at the Centre for Artificial Intelligence, Data, and Conflict,$^{[12]}$ as of May 8th, 2022, more than 3.3 million accounts were engaging with Taliban disinformation campaigns. To this day, there are expected to still be more than 120,000 accounts active on X. With X having a median engagement rate of 0.029 per post in 2024$^{[13]}$ and an average of 3 posts per account per week, we can expect the Taliban to reach roughly 0.029*120,000*3 = 10,440 accounts per week. Dividing this by seven, we find that the baseline engagement per day is 149 accounts. Since we measured NATO's engagement monthly, we'll also measure the Taliban the same way. Thus, we get the following table for the last seven months (from January 2025):
\begin{figure}[H]
    \centering
    \includegraphics[width=0.75\linewidth]{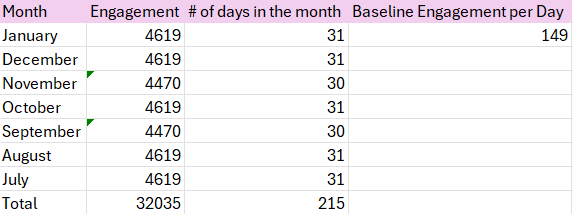}
    \caption{Expected Taliban Engagement by Month}
    \label{fig:enter-label}
\end{figure}

$\newline$Given that we can't get as much data as we would like about each individual Taliban account, this seems like a fairly good way to estimate $D_t$ with mean response.

$\newline$Putting both of our tables of mean response engagement data into one, we can finally calculate our values of $a$ for the past seven months via equation (13).
\begin{figure}[H]
    \centering
    \includegraphics[width=0.55\linewidth]{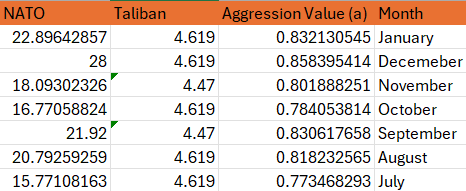}
    \caption{Aggression Table (July 2024 - January 2025)}
    \label{fig:enter-label}
\end{figure}

$\newline$Note: when we get to stochastic analysis, we will invoke the Central Limit Theorem to let $a_t$ be distributed normally with mean 0.814112 and standard deviation 0.027464. We are also given this by the fact that the composition of normally distributed variables is normally distributed. This setup isn't \textit{ideal}, but provided that NATO doesn't stop posting online and the Taliban doesn't start posting a significant amount more, we feel that this is justified. 

\subsection{Understanding the Information Environment}
While this choice may seem almost like a non-strategy, it's actually very important. Recall from the literature review that if the goodness of fit to resources, p, is greater than 1 and the second force is more efficient, and if u + v $>$ 1, then we are told to pick a = 0 until the resource market is no longer over-saturated. That is essentially the exact use case for this strategy. If for some reason NATO is overwhelmed at some point by disinformation, they should inform themselves thoroughly about how certain counter-campaigns will be received so that they don't overextend themselves - this is the intuition that follows the math!

$\newline$Note: the G7 Rapid Response Mechanism and ``Recovering Stronger" also fall under this umbrella strategy, although they admittedly take different forms.

\subsection{Mitigating Incidents}
For this approach, we want to stipulate that NATO takes one of two actions against offending disinformation campaigns:
\begin{enumerate}
    \item NATO publicly responds to refute the statements made.
    \item NATO uses its leverage to have the posts taken down and accounts banned.
\end{enumerate}

\subsubsection{The Public Recognition Strategy:}
This is possibly the most direct and aggressive tactic that NATO can take. After all, this means that everybody who follows NATO on social media will also hear the message broadcast by the insurgent group. Naturally, in this circumstance, we hold that a = 1 is a necessity. We will discuss the positives and negatives of this approach under dynamics. 

\subsubsection{The Backdoor Removal Strategy:}
This is a less direct but no less inflammatory strategy. Since making social media accounts costs little to no money, reciprocity from insurgent actors could come in the form of hundreds (if not thousands) of new social media accounts in protest. We will also hold this as an a = 1 situation, simply because of the potential for blowback. 

\section{Dynamics}
Recall our model from the literature review:
\begin{equation}
    \frac{du}{dt} = u(1-u-v)-ac_1u
\end{equation}
\begin{equation}
    \frac{dv}{dt} = pv(1-u-v)-au.
\end{equation}

$\newline$Even though we have our strategy parameter, $a$, now, we still need to find our initial condition $(u_0, v_0)$, as well as the second group's fitness to resources ($p$), and $c_i,$ the $\frac{endured}{inflicted}$ kills for population $i$.

\subsubsection{Finding p}
Before we can find $p$, we must first establish what resources are necessary to both parties, and which resources are shared. Although money is an instinctive answer, we hold that public support is the most valuable resource for both parties. Notably, the Taliban has to work much harder for it than NATO does (funny how that works). For the sake of exclusivity, we will assume that if an individual supports NATO, they do not support the Taliban, and vice versa. According to a study by the Asia Foundation$^{[14]}$ which polled Afghan citizens in 2019, only 14.9\% of respondents said that they held any sympathy for the Taliban. While we don't have a way to calculate this type of statistic for the rest of the world, we can imagine that this 14.9\% is an upper bound of the Taliban's fit to public support. Thus, we will let p = 0.149.

\subsubsection{Finding $c_1$}
To find $c_1$, we really need to qualify what we mean by a ``kill," as the civil war metaphor starts lose traction some here. If instead of making it $\frac{endured}{inflicted}$, we make it $\frac{\%detractors}{\%supporters}$, then we can reuse our statistic from 5.0.1, and let
\begin{equation}
    c_1 = \frac{0.149}{0.851} = 0.175,
\end{equation}
then we have a tractable statistic. 

\subsubsection{Finding Initial Conditions}
Finally, we just have to figure out the population density of online posts. According to a report by RadioFreeEurope / RadioFreeLiberty in 2019,$^{[15]}$ NATO paid to have 19,000 counter-information accounts in addition to it's primary methods of communication. Additionally, via Traces of Conflict,$^{[12]}$ we know that the Taliban is estimated to have around 120,000 active accounts. Thus, we will make our initial conditions
\begin{equation}
    (u_0,v_0) = (19,120).
\end{equation}
(Note that it is common to reduce orders of magnitude for simplicity. We will just interpret the values of our DE in terms of thousands of accounts.) Now we must normalize these to densities. Thus,
\begin{equation}
    (u_0, v_0) = (\frac{u_0}{u_0 + v_0}, \frac{v_0}{u_0 + v_0}) = (0.13669, 0.8633)
\end{equation}

\subsubsection{Final System of Equations}
Now that we have all of our constants and initial conditions, we can finalize the parameters in our system of equations. 
\begin{equation}
    \frac{du}{dt} = u(1-u-v) - 0.175au
\end{equation}
\begin{equation}
    \frac{dv}{dt} = 0.149v(1-u-v)-au
\end{equation}
Where a is strategy dependent, and $(u_0, v_0) = (0.13669, 0.8633)$.

\subsubsection{Dynamics of Influencer Campaigns}
Since we hold $a$ as normally distributed for influencer campaigns, we will use a = 0.814112 for this, making our system of equations:
\begin{equation}
    \frac{du}{dt} = u(1-u-v) - 0.1424696u
\end{equation}
\begin{equation}
    \frac{dv}{dt} = 0.149v(1-u-v)-0.814112u
\end{equation}

$\newline$Plotting the phase-plane, we find that NATO eventually wins since the trajectory takes the second population density to zero.  
\begin{figure}[H]
    \centering
    \includegraphics[width=0.5\linewidth]{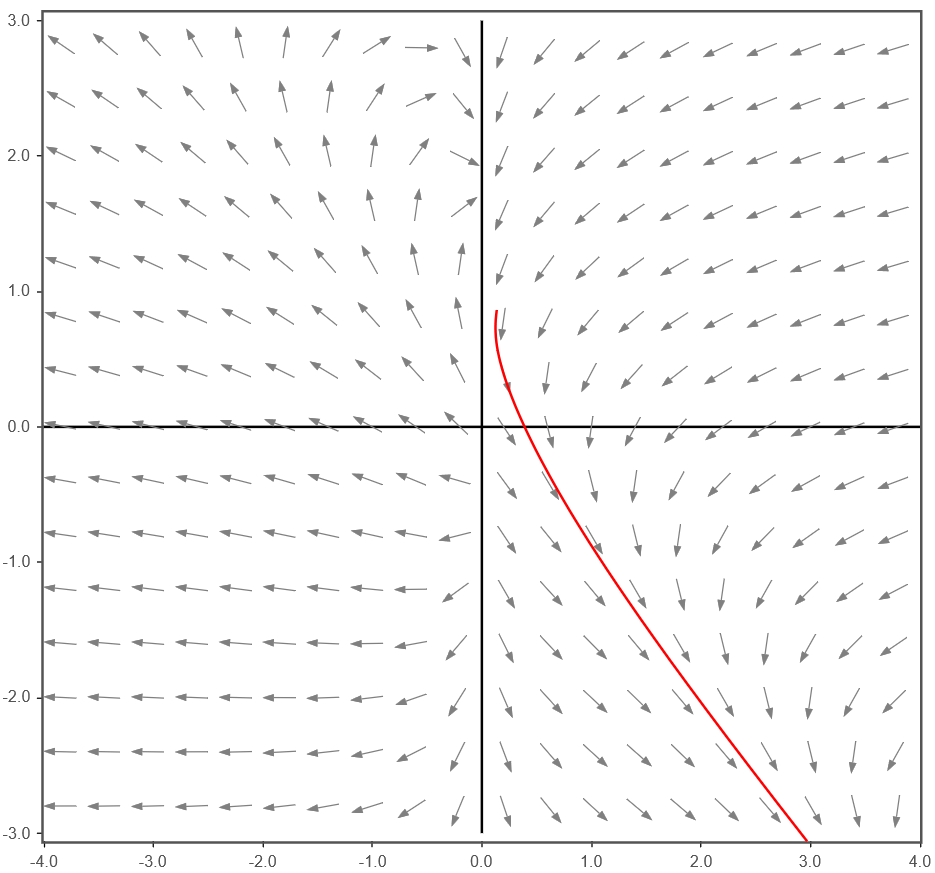}
    \caption{Phase-Plane for the Influencer Strategy}
    \label{fig:enter-label}
\end{figure}

$\newline$
According to the time plot of this trajectory, this takes about 7 ``units" of time - we'll compare this to the other strategies later.
\begin{figure}[H]
    \centering
    \includegraphics[width=0.5\linewidth]{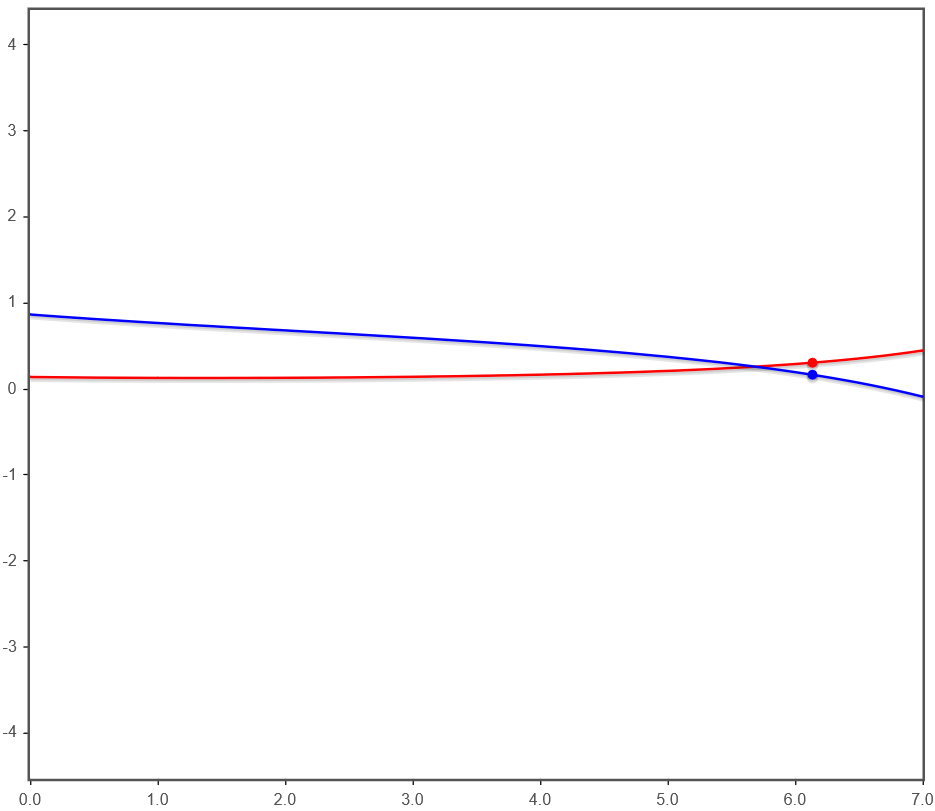}
    \caption{Time-Plot for the Influencer Strategy (Blue is the Taliban)}
    \label{fig:enter-label}
\end{figure}

\subsubsection{Dynamics of Understanding the Information Environment}
If we let a = 0, we get wildly different equations than before. Namely, 
\begin{equation}
    \frac{du}{dt} = u(1-u-v) 
\end{equation}
\begin{equation}
    \frac{dv}{dt} = 0.149v(1-u-v)
\end{equation}

$\newline$Unsurprisingly, we find that both populations stay stagnant until the other takes action. (The phase-plane is not presented because it doesn't show anything).
\begin{figure}[H]
    \centering
    \includegraphics[width=0.5\linewidth]{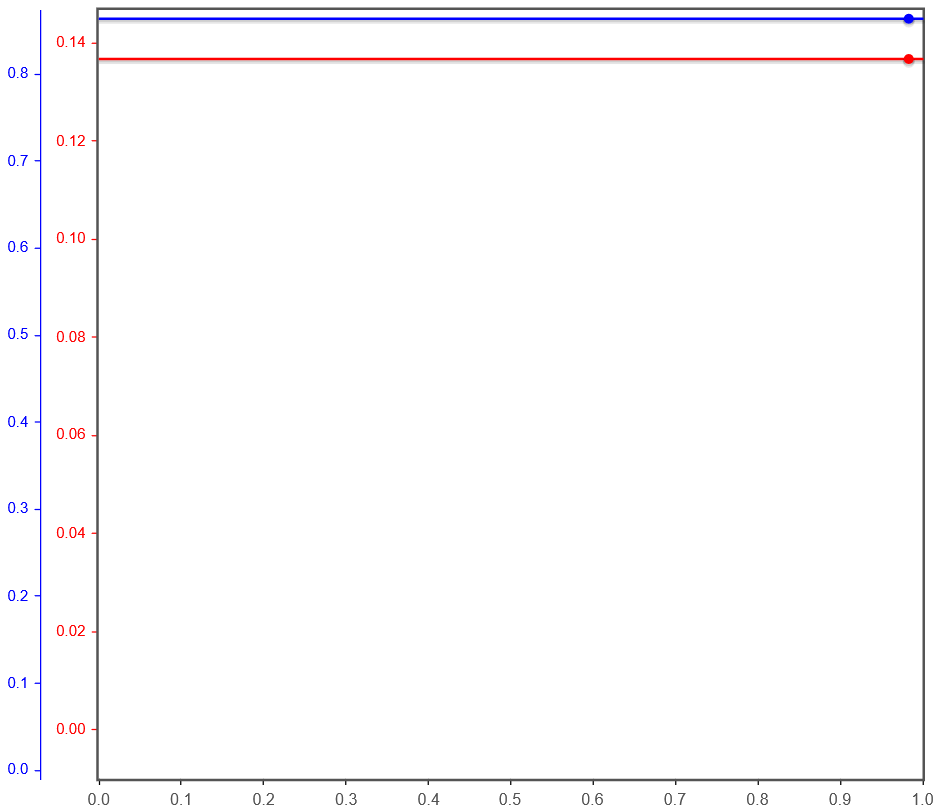}
    \caption{Time-Plot for Inaction (Blue is the Taliban)}
    \label{fig:enter-label}
\end{figure}

$\newline$The primary point of interest for this system comes from a lemma via Affili Et. Al., that if ac $\in$ (0, 1), there exists a saddle point at ($\frac{1-ac}{1+pc}pc$, $\frac{1-ac}{1+pc}$). Plugging these values in, we see that there is a saddle (roughly) at (0.020965, 0.804). That said, the implication of this is wholly dependent on the set of initial conditions that we are given from empirical data.

\subsubsection{Mitigating Incidents}
If we let a = 1, we get the following equations:
\begin{equation}
    \frac{du}{dt} = u(1-u-v) - 0.175u
\end{equation}
\begin{equation}
    \frac{dv}{dt} = 0.149v(1-u-v)-u
\end{equation}
$\newline$
Again, unsurprisingly, we find that our phase plane is just an exacerbated version of the influencer strategy:
\begin{figure}[H]
    \centering
    \includegraphics[width=0.5\linewidth]{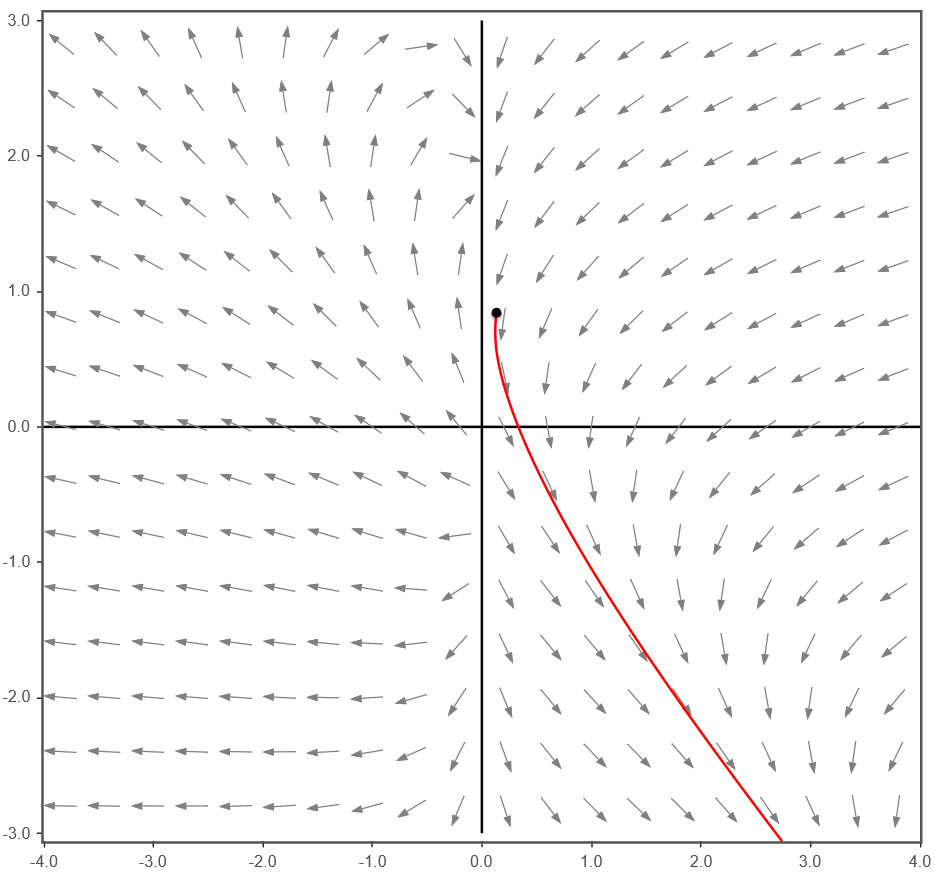}
    \caption{Phase-Plane for the Mitigation Strategy}
    \label{fig:enter-label}
\end{figure}

$\newline$
As we examine the time-plot, recall that the influencer campaign took about seven units of time to eliminate the Taliban's population density. This strategy gets it done in roughly five units. 

\begin{figure}[H]
    \centering
    \includegraphics[width=0.5\linewidth]{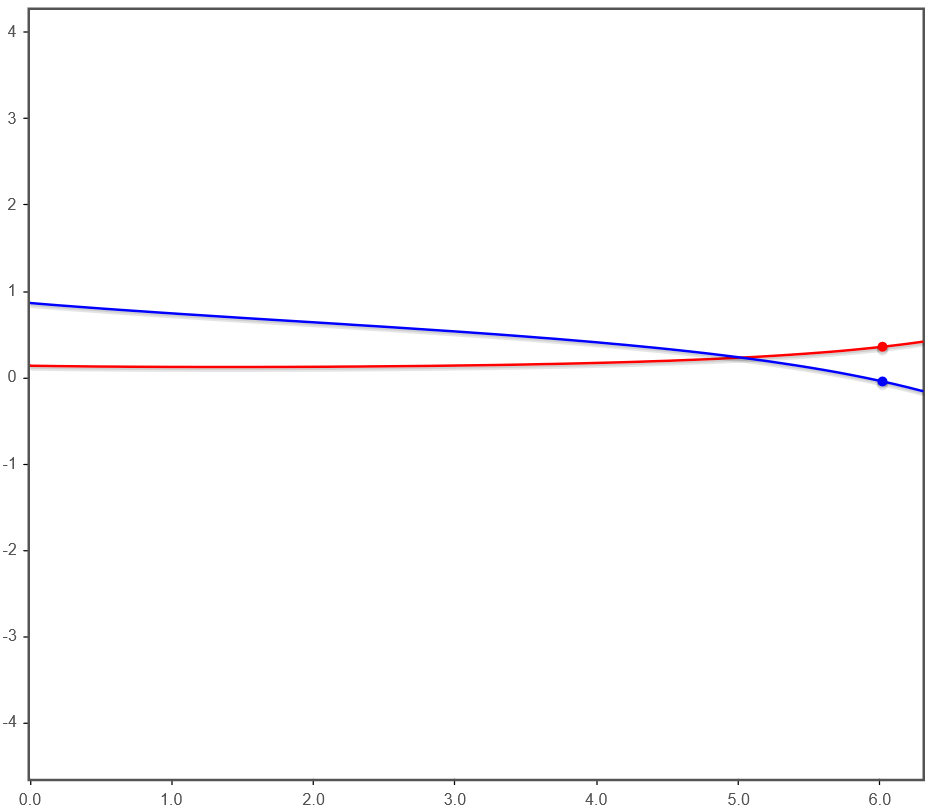}
    \caption{Time-Plot for the Mitigation Strategy (Blue is the Taliban)}
    \label{fig:enter-label}
\end{figure}

$\newline$At this point, we ask ``have we modeled everything at play?" After all, if mitigation is this effective, it only seems that the right strategy would be a purely offensive one. The truth is that we are working with an incomplete model and a lot of assumptions. The model supports a purely offensive strategy, and surely one can find historical evidence to back such policy. The model does not directly account for various types of blowback, and certainly does not account for the persistence of religious fanaticism (especially since no physical harm is being done by any strategy). Therefore, while we conclude that more aggression is preferred over inaction, we caution the reader that there may be biases of omission at play.

\section{Considering Perturbations}
In this section, we want to examine what would happen if, for some reason, either population density jumped drastically. There could be many such reasons from political dissatisfaction, to increasing leanings towards extreme beliefs, to (God forbid) increases in political violence. Regardless of why, let's examine a few of these various jumps. (Note, we omit the mitigation strategy due to its similarity to the influencer strategy.)

\subsection{A Large Taliban Spike}
Our previous initial condition was (0.13669, 0.8633), which already has the Taliban taking up a sizable proportion of the population. If we increase their population to 0.98, let's see what happens. 
\subsubsection{Influencer Campaign}
Using the same equations from (5.0.5), we get the following phase plane:
\begin{figure}[H]
    \centering
    \includegraphics[width=0.5\linewidth]{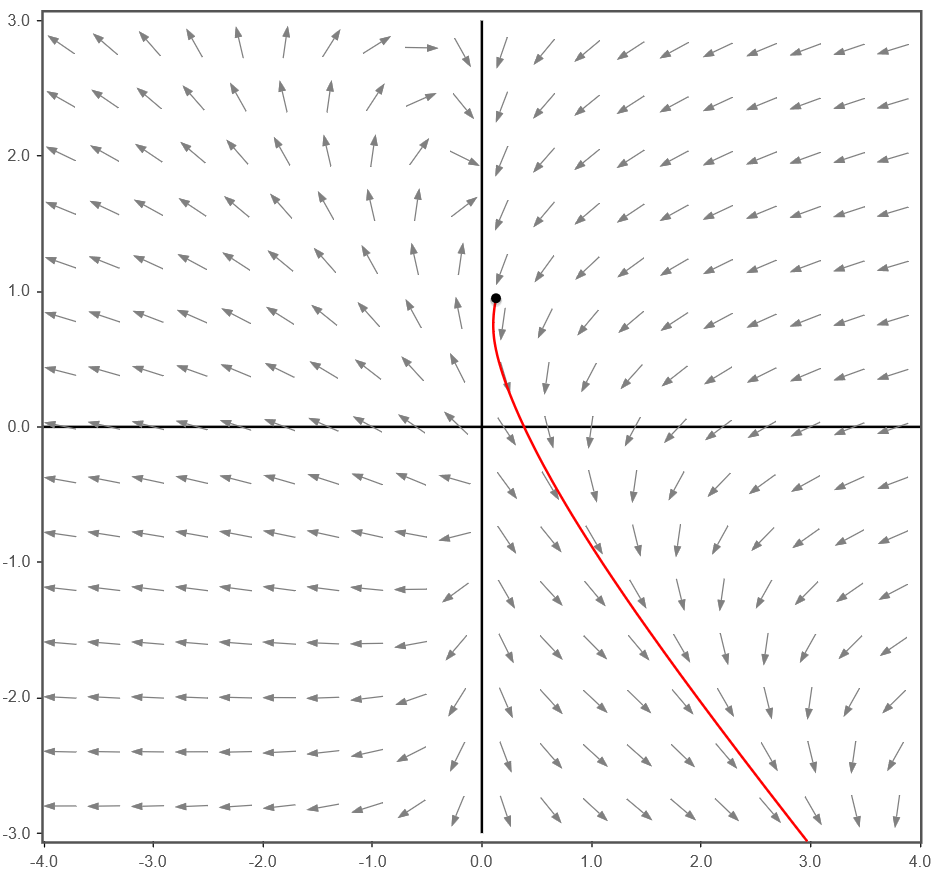}
    \caption{Phase-Plane for the Influencer Strategy}
    \label{fig:enter-label}
\end{figure}

$\newline$Notably, we don't find anything different from before except that the NATO population density gets a bit closer to zero than it did before (visible in the bowing of the trajectory). Additionally, if we examine the time-plot, we find that it merely takes a few extra units of time to eradicate the Taliban's population density - specifically, nine units total. 

\subsubsection{Understanding the Information Environment}
Using the equations from (5.0.6), we get the following time-plot:
\begin{figure}[H]
    \centering
    \includegraphics[width=0.5\linewidth]{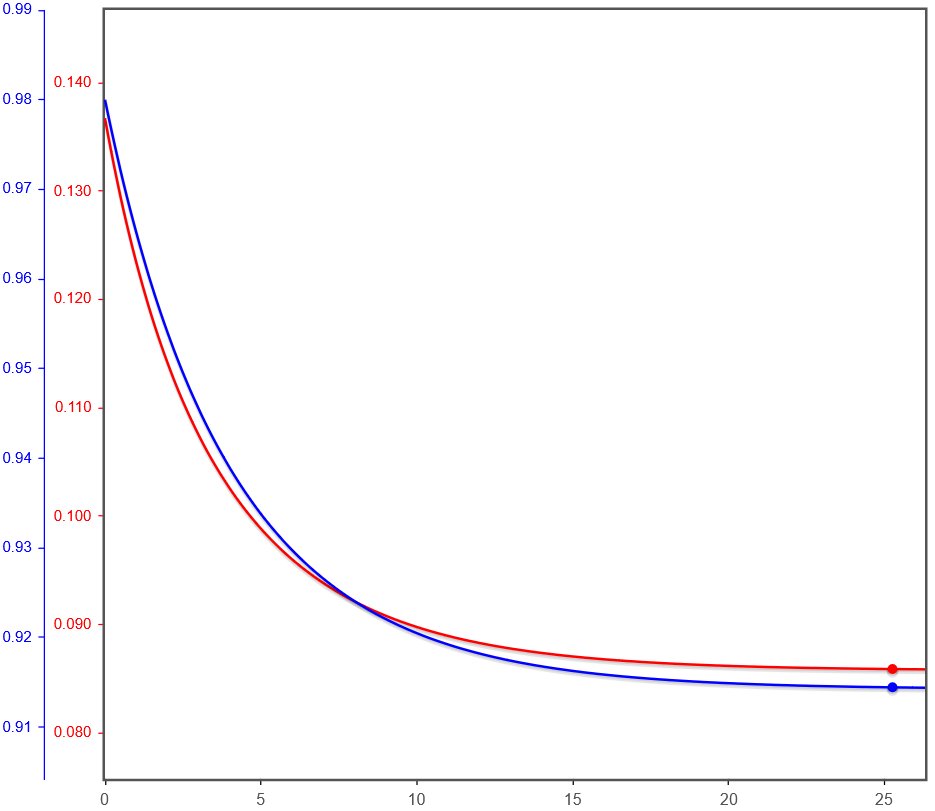}
    \caption{Time-Plot for Inaction (Blue is the Taliban)}
    \label{fig:enter-label}
\end{figure}

$\newline$
Notably, we find that starting with a higher population density of Taliban accounts causes inaction to lead to short-term losses for both parties, but still ends in the same cycle of coexistence that we found earlier.

\subsection{Moderate Taliban Spike}
This time, let's consider the halfway point between 0.98 and 0.86, $v_0$ = 0.92. 

\subsubsection{Influencer Campaign}
Once again, we find the same pattern cropping up in the phase-plane and the time-plot, again with the only difference being the amount of time it takes to happen.
\begin{figure}[H]
    \centering
    \includegraphics[width=0.5\linewidth]{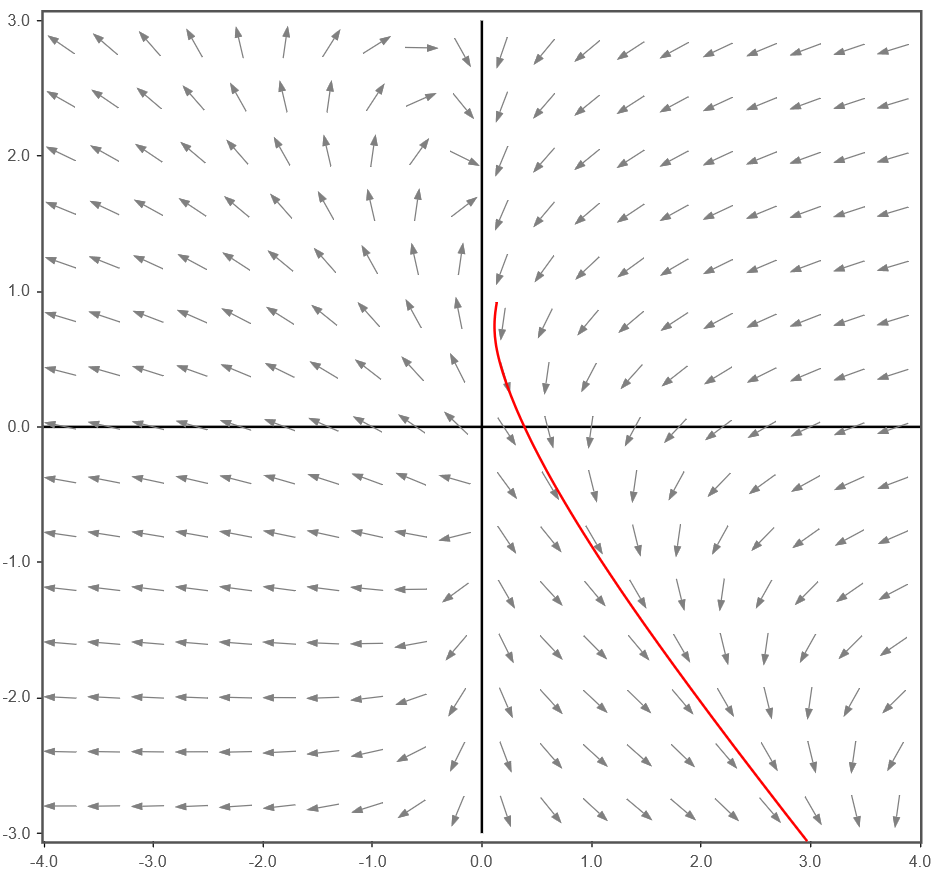}
    \caption{Phase-Plane for the Influencer Strategy}
    \label{fig:enter-label}
\end{figure}

\subsubsection{Understanding the Information Environment}
Just like last time, we find that the two trajectories fall with each other, but ultimately settle into a pattern of coexistence - only this time, they are even closer together. 
\begin{figure}[H]
    \centering
    \includegraphics[width=0.5\linewidth]{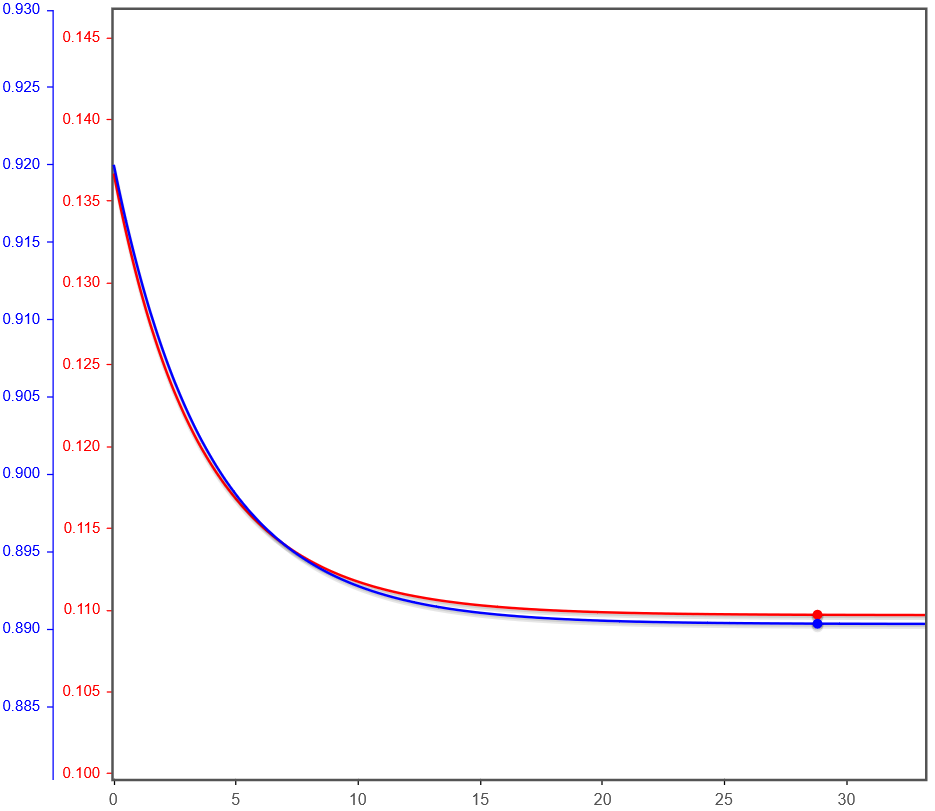}
    \caption{Time-Plot for Inaction (Blue is the Taliban)}
    \label{fig:enter-label}
\end{figure}

\subsection{Large NATO Spike}
This whole time, NATO's initial condition has been the minority of our sample space at just 0.13669. What if we boosted it to $u_0$ = 0.9?

\subsubsection{Influencer Campaign}
On the phase-plane, we see the same pattern that we have seen this whole time. However, on the time-plot, we see a slightly new trend:
\begin{figure}[H]
    \centering
    \includegraphics[width=0.5\linewidth]{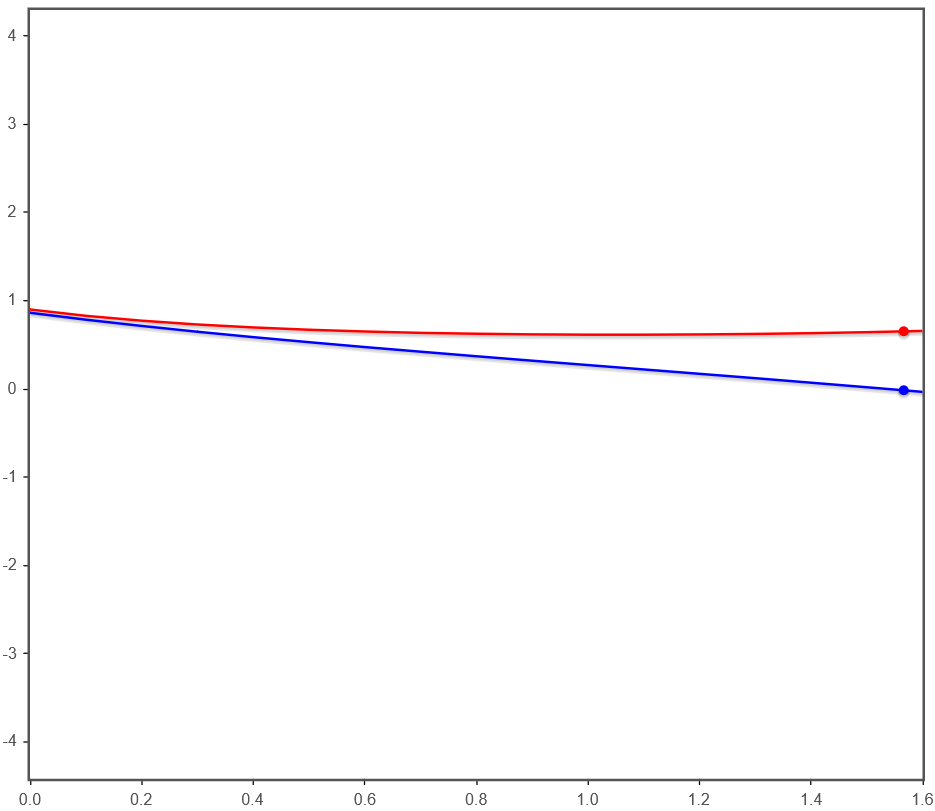}
    \caption{Time-Plot for the Influencer Strategy (Blue is the Taliban)}
    \label{fig:enter-label}
\end{figure}
$\newline$Notably, the Taliban's population density goes to zero in 1.5 units of time. 

\subsubsection{Understanding the Information Environment}
For this strategy, we find something really unique in the time-plot:
\begin{figure}[H]
    \centering
    \includegraphics[width=0.5\linewidth]{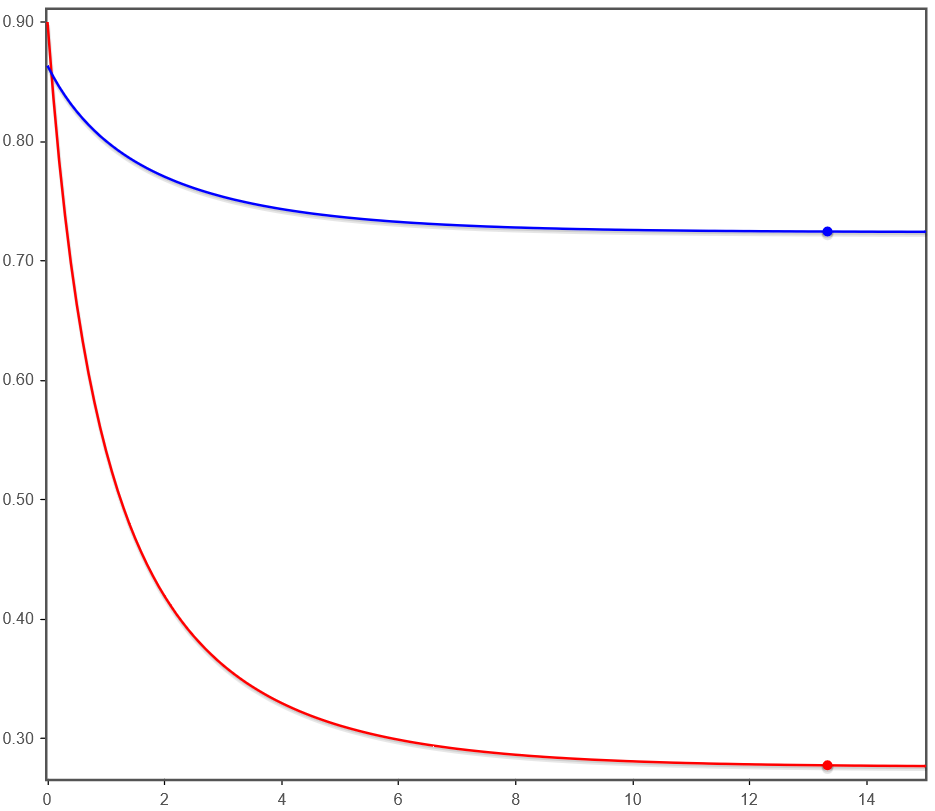}
    \caption{Time-Plot for Inaction (Blue is the Taliban)}
    \label{fig:enter-label}
\end{figure}

$\newline$Even though NATO starts off with a significant advantage (compared to the actual initial conditions), they not only lose a significant amount more of their population density than the Taliban, they also almost completely lose their population density altogether. This is a truly unexpected result with a resounding conclusion: if for some reason NATO expands their online anti-disinformation campaigns drastically, they shouldn't do it casually - otherwise, they'll lose their investment. 

\subsection{Moderate NATO Spike}
We tested out values of $u_0$ ranging between 0.13669 and 0.9, and there were not enough compelling differences between the values we've already tested and the results we've already shown to include any of the examples.  

\section{Policy Conclusions}
From the various strategies we have tested, we find that the optimal strategy for the current state-of-the-world is to be fully aggressive, even to the point of directly calling out disinformation campaigns as needed. That said, the marginal benefit between mitigation and an influencer campaign is minimal, so NATO should balance its strategy selection accordingly. One potential policy action to fit this balance is the one discussed in section (2.4) - the NATO Afghan Youth for Change program. By targeting individuals in the Afghan diaspora to become NATO influencers, NATO can supplant their own agenda into communities which might otherwise be earmarked by insurgent groups for grooming into terror cells. This is a very aggressive strategy with significantly smaller potential for blowback than direct mitigation, which makes it extremely optimal for implementation in the near future.

$\newline$
We also find that if NATO expands its campaign for countering disinformation with bot accounts, influencers, or other options, it must absolutely not do it passively - otherwise it runs the risk of atrophying its own resources away.

\vfill

\section{References}
[1] Nato. (2023, March 3). NATO’s approach to counter information threats. NATO. 

https://www.nato.int/cps/in/natohq/topics\_219728.htm\#how 

$\newline$
[2] Nato. (2025). Protect the future. NATO. 

https://www.nato.int/protect-the-future/ 

$\newline$
[3] Pew Research Center. (2024, September 17). Social Media and News Fact sheet. Pew Research Center. 

https://www.pewresearch.org/journalism/fact-sheet/social-media-and-news-fact-sheet/ 

$\newline$
[4] Hughes, K. (2019, May). Fact sheet – information environment assessment (IEA). NATO. 

https://www.act.nato.int/wp-content/uploads/2023/05/2019\_05\_IEA.pdf 

$\newline$
[5] Information environment assessment capability ... NATO. (2023, April 3). 

https://www.act.nato.int/article/information-environment-assessment-capability-programme-plan-initiated 

$\newline$
[6] Government of Canada. (2024, October 3). Government of Canada. GAC. 

https://www.international.gc.ca/transparency-transparence/rapid-response-mechanism-mecanisme-reponse-rapide/index.aspx?lang=eng 

$\newline$
[7] NATO Public Diplomacy Division. NATO. (2024, November 8). 

https://www.nato.int/structur/pdd/2025/2025-content-guidelines-nato-grants.pdf 

$\newline$
[8] Annual report 2022. G7 2022 Annual Report. (2022). 

https://www.international.gc.ca/transparency-transparence/assets/pdfs/rapid-response-mechanism-mecanisme-reponse-rapide/g7-rrm-2022-annual-report-en.pdf 

$\newline$
[9] Akram, M., Nasar, A., \& Perveen, S. (2024, July). Decoding social media’s role in Taliban 2.0 and its ... Asia Pacific Issues. 

https://www.eastwestcenter.org/sites/default/files/2024-07/API No. 168 . Taliban 2.0 \_0.pdf 

$\newline$[10] Instagram Audit Report: Phlanx Influencer directory. Instagram Audit Report | Phlanx Influencer Directory. (2025). 

https://phlanx.com/auditor/instagram/1866a455cdeb431e?export=pdf 

$\newline$[11] GeeksforGeeks. (2024b, April 30). Overview on Lotka-Volterra model of predator-prey relationship - GeeksforGeeks. 

https://www.geeksforgeeks.org/lotka-volterra-model-of-predator-prey-relationship/\#what-is-the-equation-for-the-predatorprey-model

$\newline$[12] Courchesne, L., Rasikh, B., McQuinn, B., \& Buntain, C. (2022). Powered by Twitter? The Taliban’s Takeover of Afghanistan. Traces of Conflict. 

https://www.tracesofconflict.com/\_files/ugd/17ec87\_19ecafa8cf1046af8554251bce0aaf6f.pdf 

$\newline$[13] Lauron, S. (2024, June 15). What is a good engagement rate on Twitter?. Rival IQ. 

https://www.rivaliq.com/blog/good-engagement-rate-twitter/ 

$\newline$[14] Maley, W. (2021, September 17). The public relations of the Taliban: Then and now. ICCT. 

https://icct.nl/publication/public-relations-taliban-then-and-now 

$\newline$[15] RFE/RL. (2019, December 6). NATO-linked report says social-media platforms failing on manipulation. RadioFreeEurope/RadioLiberty. 

https://www.rferl.org/a/nato-linked-report-says-social-media-platforms-failing-on-manipulation/

30310914.html 

$\newline$[16] Affili, E., Dipierro, S., Rossi, L., \& Valdinoci, E. (2024). A new Lotka-Volterra model of competition with strategic aggression: Civil wars when strategy comes into play. Birkhäuser. 

$\newline$Additional thanks to Dr. John Scherpereel, Dr. Bernd Kaussler, Sean Tarter, Vikram Lothe, Roland Allen, and Savannah Whitley for their insightful conversations and grammar edits. 

\end{document}